\documentclass[12pt]{article}

\usepackage{mathrsfs}
\usepackage{color}
\usepackage{amssymb, amsthm, amsfonts, amsxtra, amsmath}
\usepackage{mathabx}

\usepackage{latexsym}
\usepackage[all]{xy}
\usepackage{graphics}
\usepackage{latexsym}
\usepackage{makeidx}
\usepackage{rotating}

\theoremstyle{change}

\newtheorem{Theorem}{Theorem}[section]
\newtheorem{Def}[Theorem]{Definition}
\newtheorem{Lem}[Theorem]{Lemma}
\newtheorem{Prop}[Theorem]{Proposition}

\newtheorem{Def-Prop}[Theorem]{Definition-Proposition}
\newtheorem{Not}[Theorem]{Notation}

\date{}

\begin{document}

\hyphenation{Wo-ro-no-wicz}

\title{On the construction of quantum homogeneous spaces from $^*$-Galois objects}
\author{Kenny De Commer\footnote{Supported in part by the ERC Advanced Grant 227458
OACFT ``Operator Algebras and Conformal Field Theory" }\\ \small Department of mathematics,  University of Tor Vergata, Rome\\
\small Via della Ricerca Scientifica 1, 00133 Roma, Italy\\ \\ \small e-mail: decommer@mat.uniroma2.it}
\maketitle

\newcommand{\acnabla}{\nabla\!\!\!{^\shortmid}}
\newcommand{\undersetmin}[2]{{#1}\underset{\textrm{min}}{\otimes}{#2}}
\newcommand{\otimesud}[2]{\overset{#2}{\underset{#1}{\otimes}}}
\newcommand{\qbin}[2]{\left[ \begin{array}{c} #1 \\ #2 \end{array}\right]_{q^2}}

\newcommand{\otimesmin}{\underset{\textrm{min}}{\otimes}}
\newcommand{\bigback}{\!\!\!\!\!\!\!\!\!\!\!\!\!\!\!\!\!\!\!\!\!\!\!\!}

\abstract{\noindent In this note we construct bi-$^*$-Galois objects linking the quantized universal enveloping algebras associated to the Lie groups $SU(2)$, $E(2)$ and $SU(1,1)$, where $E(2)$ denotes the Lie group of Euclidian transformations of the plane, and we show how one can create (formal) quantum homogeneous spaces for these quantum groups by integrating the associated Miyashita-Ulbrich action on certain subquotient $^*$-algebras.}\vspace{0.5cm}


\noindent AMS 2000 \emph{Mathematics subject classification}: 17B37; 16W30; 81R50; 46L65


\section*{Introduction}

\noindent This is part of a series of papers (\cite{DeC1}) devoted to an intriguing correspondence between the quantizations of $SU(2)$, $\tilde{E}(2)$ and $SU(1,1)$, where $\tilde{E}(2)$ denotes the non-trivial double cover of the Lie group of Euclidian transformations of the plane. In a sense, their duals form a trinity of `Morita equivalent locally compact quantum groups'. There then exists a `linking quantum groupoid' combining these three quantum groups into one global structure, and it is important to understand for example the (co)representation theory of this object.\\

\noindent In this paper, we will study an infinitesimal structure associated to this linking quantum groupoid. This discussion will only be preliminary, in the sense that we investigate the elementary algebraic structure, and do not consider the more delicate issues (concerning for example the spectrum of Casimir elements in `admissible' Hilbert space representations), which will be treated elsewhere. Indeed, the main goal of this article is simply to collect some basic results concerning this structure, some of which were already observed in the literature, and to show how its \emph{elementary} representation theory can be used to give a novel construction of quantum homogeneous spaces for each of the constituent quantum groups.\\

\noindent In the \emph{first section}, we will start by introducing the theory of bi-$^*$-Galois objects between $U_q(su(2))$, $U_q(e(2))$ and $U_q(su(1,1))$. In fact, the bi-Galois objects themselves are already well-known in the literature (their first appearance seems to be in \cite{Gun1}, while in \cite{Kas1} higher-dimensional generalizations are considered (see also \cite{Did1}, \cite{Aub1} and \cite{Mas1})). The only new observation then is that these bi-Galois objects allow for a compatible $^*$-structure (in case $q$ is real). While this observation in itself is quite trivial, it has, as one can expect, most important consequences if one considers representation theoretic issues. It also explains the appearance of the third QUE algebra $U_q(su(1,1))$, which was not present in these earlier papers, since, when neglecting the $^*$-structure, this is just an isomorphic copy of $U_q(su(2))$. Apart from this, we introduce an auxiliary organizing structure, called a \emph{co-linking weak Hopf $^*$-algebra}, which allows us to treat the various bi-$^*$-Galois objects in a unified way. This structure is a small adaptation of the notion of a total Hopf Galois system (\cite{Gru1}), which was itself an enhancement of the notion of a Hopf Galois system, introduced in \cite{Bic1}. The further improvement is that we can simplify the axiom system of such a total Hopf Galois system by using the language of weak Hopf algebras (\cite{Boh1}). In practice however, one always uses the same techniques for any of these notions.\\


\noindent In the \emph{second section}, we introduce a natural notion of \emph{quantum Casimir element} inside these bi-$^*$-Galois objects, which will be a certain self-adjoint element in the center. We then consider the quotient $^*$-algebras, obtained by evaluating the Casimir element at a particular real value. These quotients carry a natural module $^*$-algebra structure, obtained from the Miyashita-Ulbrich action on the original bi-$^*$-Galois objects. Finally, considering sub-$^*$-algebras, we obtain our quantum homogeneous spaces equipped with their natural infinitesimal action. We also determine which of these quantum homogeneous spaces arise from $^*$-coideals.\\

\noindent \qquad \emph{Remarks on notation and conventions}\\

\noindent For the rest of the paper, we fix a real number $0<q<1$. We then denote \[\lambda = (q-q^{-1})^{-1} < 0.\]

\noindent By $\iota$ we always mean the identity map, and by $\otimes$ the tensor product between vector spaces over $\mathbb{C}$.\\

\noindent We will use the Sweedler notation for comultiplications (see \cite{Swe1}). Also, we will only work with Hopf algebras over $\mathbb{C}$ which have invertible antipodes.\\

\noindent We will need a lot of structures which are defined very similarly to each other. Then the names for the structures are often indexed, and when multiple structures are used together, we will index the objects associated to these structures with the corresponding index. However, when the structures appear isolated, we will refrain from indexing any of its associated structure. Also, when we index something with two indices which are the same, we will sometimes take the liberty of indexing with just one times this index symbol. We hope that in practice, these conventions will cause no confusion. \\

\section{Bi-$^*$-Galois objects for $U_q(su(2))$, $U_q(e(2))$ and $U_q(su(1,1))$.}

\begin{Def} Let $\mu,\nu\in \mathbb{R}$. The $^*$-algebra $U_q(\mu,\nu)$ is defined as the universal algebra generated by four elements $K,K^{-1},E$ and $F$ satisfying the commutation relations \[K^{-1}K = 1 =  KK^{-1}, \;\;\;\; KE = qEK,\;\;\;\;  KF = q^{-1}FK \] and \[\lbrack E,F\rbrack = \lambda(\mu K^{2} - \nu K^{-2}),\] endowed with the $^*$-structure determined by $K^*=K$ and $E^*=F$.\\

\noindent When $\mu=\nu$, the $^*$-algebra $U_q(\mu) = U_q(\mu,\mu)$ can be turned into a Hopf $^*$-algebra, with the comultiplication determined by
\[ \Delta(E) = E\otimes K+K^{-1}\otimes E\] and \[\Delta(K) = K\otimes K.\]

\end{Def}

\noindent It is not difficult to see that, by rescaling $E$, $F$ and $K$ by some positive number, all Hopf $^*$-algebras $U_q(\mu)$ with $\mu>0$ are isomorphic to the Hopf $^*$-algebra $U_q(1)=U_q(su(2))$ (\cite{Kli1}), which we will then also denote as $U_q(+)$. Similarly, all $U_q(\mu)$ with $\mu<0$ are isomorphic to the Hopf $^*$-algebra $U_q(-1)=U_q(su(1,1))$ (\cite{Kli1}), which we will then also denote as $U_q(-)$. The remaining Hopf $^*$-algebra $U_q(0)$ can be interpreted as the quantum universal enveloping algebra of the Lie algebra $e(2)$ associated to the Lie group of Euclidian transformations of the plane (\cite{Vak1}). We note that in some contexts, it is better not to rescale the parameter $\mu$. For example, with regard to the contraction procedure (\cite{Cel1}), one considers $U_q(0)$ as $\lim_{\mu\rightarrow 0} U_q(\mu)$.\\

\noindent When $\mu\neq \nu$, the $^*$-algebra $U_q(\mu,\nu)$ does not always possess a good Hopf $^*$-algebra structure. However, it \emph{does} possess a natural $U_q(\mu)$-$U_q(\nu)$-bicomodule $^*$-algebra structure, which is in some sense the smallest possible such structure which is well enough behaved to allow one to go back and forth between $U_q(\nu)$ and $U_q(\nu)$ without losing information. In less vague terms, the $U_q(\mu,\nu)$ are bi-$^*$-Galois objects, whose definition we now recall.

\begin{Def} (\cite{Sch1}, Definition 3.4; \cite{Bic0}, Definition 4.1.1) Let $(H_0,\Delta_0)$ and $(H_1,\Delta_1)$ be two Hopf $^*$-algebras. A \emph{bi-$^*$-Galois object} $(H_{01},\Delta_{01}^0,\Delta_{01}^1)$ between $(H_0,\Delta_0)$ and $(H_1,\Delta_1)$ consists of a non-trivial unital $^*$-algebra $H_{01}$ (i.e. $1\neq 0$) and \begin{enumerate} \item[1.] a left $H_{0}$-comodule $^*$-algebra structure (= \emph{left coaction}) \[\Delta_{01}^0:H_{01}\rightarrow H_0\otimes H_{01},\] and \item[2.] a right $H_{1}$-comodule $^*$-algebra structure (= \emph{right coaction}) \[\Delta_{01}^1:H_{01}\rightarrow H_{01}\otimes H_1,\]\end{enumerate} such that \begin{enumerate}
\item[3.] the \emph{left canonical map} \[H_{01}\otimes H_{01}\rightarrow H_{01}\otimes H_{1}: x\otimes y\rightarrow (x\otimes 1)\Delta_{01}^1(y)\] is bijective,
\item[4.] the \emph{right canonical map} \[H_{01}\otimes H_{01}\rightarrow H_0\otimes H_{01}: x\otimes y \rightarrow \Delta_{01}^0(x)(1\otimes y)\] is a bijection,\end{enumerate} and the maps $\Delta_{01}^0$ and $\Delta_{01}^1$ commute: \[(\iota\otimes \Delta_{01}^1)\Delta_{01}^0 = (\Delta_{01}^0\otimes \iota)\Delta_{01}^1.\]

\noindent When $(H_1,\Delta_1)$ (resp.~ $(H_0,\Delta_0)$) is a Hopf $^*$-algebra, and $(H_{01},\Delta_{01}^1)$ (resp.~ $(H_{01},\Delta_{01}^0)$) satisfies the even (resp.~ uneven) conditions (and disregarding the final commutativity statement), it is called a \emph{right (resp.~ left) $^*$-Galois object} for $(H_1,\Delta_1)$ (resp. $(H_0,\Delta_0)$).

\end{Def}

\noindent One of the main results of \cite{Sch1} states, in the non-$^*$-case, that the existence of a bi-Galois object between two Hopf algebras coincides with the existence of a \emph{monoidal equivalence} between their comodule categories (see \cite{Bic0} for the $^*$-case). Given this categorical interpretation, one sees that such Hopf $^*$-algebras are still very much the same, and can in some sense be interpreted as different algebraic implementations (or representations) of `the same underlying quantum group'. The \emph{other} main result of \cite{Sch1}, in its content closely connected with the Tannaka-Krein reconstruction theorems, is that any right Galois object for a Hopf algebra \emph{can be completed to a bi-Galois object}. In particular, from a right Galois object $(H_{01},\Delta_{01}^1)$ for $(H_1,\Delta_1)$ one constructs a \emph{new} Hopf algebra $(H_0,\Delta_0)$, which in many cases of interest is indeed non-isomorphic to the original one. Finally, a bi-Galois object is even completely determined by its associated right Galois object (up to isomorphism of $(H_0,\Delta_0)$ of course).\\

\noindent The following lemma will allow us to build a bi-$^*$-Galois object structure on the $U_q(\mu,\nu)$.

\begin{Lem}\label{LemCoas} For $\mu,\upsilon,\nu\in \mathbb{R}$, there exists a unique unital $^*$-homomorphism \[\Delta_{\mu\nu}^{\upsilon}: U_q(\mu,\nu)\rightarrow U_q(\mu,\upsilon)\otimes U_q(\upsilon,\nu)\]such that \[\Delta_{\mu\nu}^{\upsilon}(E_{\mu\nu}) = E_{\mu\upsilon}\otimes K_{\upsilon\nu}+K_{\mu\upsilon}^{-1}\otimes E_{\upsilon\nu}\] and \[\Delta_{\mu\nu}^{\upsilon}(K_{\mu\nu}) = K_{\mu\upsilon}\otimes K_{\upsilon\nu}.\] Then these $^*$-homomorphisms satisfy the generalized coassociativity condition \[(\Delta_{\mu\upsilon}^{\omega}\otimes \iota)\Delta_{\mu\nu}^{\upsilon} = (\iota\otimes \Delta_{\omega\nu}^{\upsilon})\Delta_{\mu\nu}^{\omega}.\]  \end{Lem}

\begin{proof} To see if the $\Delta$-maps are well-defined, we should only see if the $\Delta_{\mu\nu}^{\upsilon}$ respect the commutation relations, since it is immediately seen that these maps will then preserve the $^*$-structure.\\

\noindent We only make the computation for the final commutation relation between $E_{\mu\nu}$ and $F_{\mu\nu}$, since the preservation of the $q$-commutation relations is apparent on sight:

\begin{eqnarray*} \lbrack \Delta_{\mu\nu}^{\upsilon}(E_{\mu\nu}),\Delta_{\mu\nu}^{\upsilon}(F_{\mu\nu})\rbrack &=& \quad\! (E_{\mu\upsilon}F_{\mu\upsilon}\otimes K_{\upsilon\nu}K_{\upsilon\nu}
- F_{\mu\upsilon} E_{\mu\upsilon}\otimes K_{\upsilon\nu}K_{\upsilon\nu})\\ &&+ (K_{\mu\upsilon}^{-1}F_{\mu\upsilon}\otimes E_{\upsilon\nu} K_{\upsilon\nu}- F_{\mu\upsilon}K_{\mu\upsilon}^{-1}\otimes K_{\upsilon\nu}E_{\upsilon\nu}) \\&&+ (E_{\mu\upsilon}K_{\mu\upsilon}^{-1}\otimes K_{\upsilon\nu}F_{\upsilon\nu}-K_{\mu\upsilon}^{-1} E_{\mu\upsilon}\otimes F_{\upsilon\nu}K_{\upsilon\nu})\\ &&
+ (K_{\mu\upsilon}^{-1} K_{\mu\upsilon}^{-1}\otimes E_{\upsilon\nu}F_{\upsilon\nu}- K_{\mu\upsilon}^{-1}K_{\mu\upsilon}^{-1} \otimes F_{\upsilon\nu}E_{\upsilon\nu})\\ \\&=& \quad \! \lbrack E_{\mu\upsilon},F_{\mu\upsilon}\rbrack\otimes K_{\upsilon\nu}^2 \\ && + (K_{\mu\upsilon}^{-1}F_{\mu\upsilon}\otimes E_{\upsilon\nu} K_{\upsilon\nu}- q^{-1}\cdot q\cdot  K_{\mu\upsilon}^{-1} F_{\mu\upsilon}\otimes E_{\upsilon\nu} K_{\upsilon\nu}) \\ && + (E_{\mu\upsilon}K_{\mu\upsilon}^{-1}\otimes K_{\upsilon\nu}F_{\upsilon\nu}- q^{-1}\cdot q \cdot  E_{\mu\upsilon}K_{\mu\upsilon}^{-1} \otimes K_{\upsilon\nu}F_{\upsilon\nu}) \\ &&+ K_{\mu\upsilon}^{-2}\otimes \lbrack E_{\upsilon\nu},F_{\upsilon\nu}\rbrack\\ \\ &=& \lambda (\mu K_{\mu\upsilon}^2 - \upsilon K_{\upsilon\nu}^{-2})\otimes K_{\upsilon\nu}^2 + K_{\mu\upsilon}^{-2}\otimes \lambda (\upsilon K_{\upsilon\nu}^2 - \nu K_{\upsilon\nu}^{-2}) \\ &=&\lambda \mu K_{\mu\upsilon}^2\otimes K_{\upsilon\nu}^2- \lambda \nu K_{\mu\upsilon}^{-2}\otimes K_{\upsilon\nu}^{-2}\\ &=& \lambda (\mu \Delta_{\mu\nu}^{\upsilon}(K_{\mu\nu}^2) - \nu\Delta_{\mu\nu}^{\upsilon}(K_{\mu\nu}^{-2})) \\ &=& \Delta_{\mu\nu}^{\upsilon}(\lbrack E_{\mu\nu},F_{\mu\nu}\rbrack). \end{eqnarray*}

\noindent The generalized coassociativity condition should then only be checked on generating elements. But, if we forget the index, then we see that we are just looking at one particular and well-known 4-dimensional coalgebra. This concludes the proof.

\end{proof}

\noindent We will then, for $x\in U_q(\mu,\nu)$, denote \[\Delta_{\mu\nu}^{\upsilon}(x) = x_{(1)\mu\upsilon}\otimes x_{(2)\upsilon\nu} \in U_q(\mu,\upsilon)\otimes U_q(\upsilon,\nu).\]

\begin{Lem} For all $\mu,\nu\in \mathbb{R}$, there exists a bijective anti-homomorphism \[S_{\mu\nu}:U_q(\mu,\nu)\rightarrow U_q(\nu,\mu)\] such that \begin{eqnarray*} S_{\mu\nu}(E_{\mu\nu}) &=& -qE_{\nu\mu},\\ S_{\mu\nu}(F_{\mu\nu})&=&-q^{-1}F_{\nu\mu}, \\ S_{\mu\nu}(K_{\mu\nu}) &=& K_{\nu\mu}^{-1}.\end{eqnarray*}
\end{Lem}

\begin{proof} Let us check again that the commutation relation between $E$ and $F$ is preversed, the other statements being verifiable on sight. We compute: \begin{eqnarray*} S_{\mu\nu}(\lbrack E_{\mu\nu},F_{\mu\nu}\rbrack) &:=&  \lbrack S_{\mu\nu}(F_{\mu\nu}),S_{\mu\nu}(E_{\mu\nu})\rbrack \\ &=& - \lbrack E_{\nu\mu},F_{\nu\mu}\rbrack \\ &=& -\lambda (\nu K_{\nu\mu}^2-\mu K_{\nu\mu}^{-2}) \\ &=& \lambda (\mu S_{\mu\nu}(K_{\mu\nu}^2) -\nu S_{\mu\nu}(K_{\mu\nu}^{-2})) \\ &=& S_{\mu\nu}(\lambda (\mu K_{\mu\nu}^2- \nu K_{\mu\nu}^{-2})).\end{eqnarray*}

\end{proof}

\begin{Lem}\label{LemAnt} For all $\mu,\nu\in \mathbb{R}$ and $x\in U_q(\mu,\mu)$, we have that \[S_{\mu\nu}(x_{(1)\mu\nu})x_{(2)\nu\mu} = \varepsilon_{\mu}(x)1_{\nu\mu}\qquad \textrm{and}\qquad  x_{(1)\mu\nu}S_{\nu\mu}(x_{(2)\nu\mu}) = \varepsilon_{\mu}(x)1_{\mu\nu}.\]\end{Lem}

\begin{proof} We only have to see if the above identities hold true for $x$ a generator. For $x=K,K^{-1}$, the identities are immediate. For $x=E_{\mu}$, we have \begin{eqnarray*} S_{\mu\nu}(E_{\mu(1)\mu\nu})E_{\mu(2)\nu\mu} &=& S_{\mu\nu}(E_{\mu\nu})\cdot K_{\nu\mu} + S_{\mu\nu}(K_{\mu\nu}^{-1})\cdot E_{\nu\mu} \\ &=& -qE_{\nu\mu}K_{\nu\mu} + K_{\nu\mu}E_{\nu\mu} \\ &=& -q E_{\nu\mu}K_{\nu\mu} +q E_{\nu\mu}K_{\nu\mu} \\ &=& 0 \\ &=&  \varepsilon_{\mu}(E_{\mu})1_{\nu\mu}.\end{eqnarray*} The other equalities to check are all similar.\end{proof}

\noindent We can now easily prove the following proposition.

\begin{Prop} Every $U_q(\mu,\nu)$ is a $U_q(\mu)$-$U_q(\nu)$-bi-$^*$-Galois object.
\end{Prop}
\begin{proof} Using the previous lemmas, one immediately verifies that the $U_q(\mu,\nu)$, where $\mu,\nu$ range over two specific values, satisfy the definitions of a Hopf-Galois system as in Definition 1.1 of \cite{Bic1}. The only thing which still needs clarification is that the $U_q(\mu,\nu)$ are not zero for $\mu\neq \nu$, but this will follow from the proof of Proposition \ref{PropUniv}. The Proposition now follows from Theorem 1.2 of \cite{Bic1}.
\end{proof}

\noindent In fact, we prefer to use the language of \emph{co-linking weak Hopf $^*$-algebras}, introduced in the next definition, as opposed to that of the Hopf-Galois systems of \cite{Bic1}, since the latter are not very symmetrical in their definition (a deficit well acknowledged by the author of that paper, whose only goal was to simply capture the essence of the structure to create bi-Galois objects). These co-linking weak Hopf-$^*$-algebras can be seen as specializations of Takeuchi's pre-equivalences (or strict Morita-Takeuchi-contexts as they are now called). The notion of a co-linking weak Hopf $^*$-algebra can also be shown to be equivalent with that of a \emph{total Hopf-Galois system} of \cite{Gru1} (equipped with a $^*$-structure), but using the language of weak Hopf algebras makes the definition somewhat more concise. The proof of the equivalence between these two concepts is essentially the one of Proposition \ref{PropWHA}.\\

\begin{Def}\label{DefWHo} For $i,j\in\{0,1\}$, let the $H_{ij}$ be four non-trivial $^*$-algebras, equipped with eight unital $^*$-homomorphisms $\Delta_{ij}^k: H_{ij}\rightarrow H_{ik}\otimes H_{kj}$, $i,j,k\in\{0,1\}$. We then call this collection a \emph{co-linking weak Hopf $^*$-algebra} when $H := \underset{i,j\in \{0,1\}}{\oplus} H_{ij}$, the direct sum $^*$-algebra, together with the map \[\Delta_H: H\rightarrow H\otimes H:\]\[H\supseteq H_{ij}\ni x\rightarrow \Delta_{ij}^0(x)+\Delta_{ij}^1(x) \in H_{i0}\otimes H_{0j}+H_{i1}\otimes H_{1j}\subseteq H\otimes H,\] forms a weak Hopf $^*$-algebra (\cite{Boh1}, Definition 2.1 and Definition 4.1, disregarding the condition of finite dimensionality). \end{Def}

\noindent One may thus interpret such a co-linking weak Hopf $^*$-algebra as the (algebraic) function space on some groupoid with two objects (indeed, in \cite{Bic1} it was already observed that Hopf-Galois systems form incomplete descriptions of such objects).\\

\noindent One can further show that a bi-$^*$-Galois object between two Hopf $^*$-algebras can be completed in an essentially unique way to a co-linking weak Hopf $^*$-algebra between the two Hopf $^*$-algebras, and that, conversely, the $(H_{01},\Delta_{01}^0,\Delta_{01}^1)$-part of a co-linking weak Hopf $^*$-algebra gives a $(H_0,\Delta_{00}^0)$-$(H_{1},\Delta_{11}^1)$-bi-$^*$-Galois object (using for example Theorem 6.1 of \cite{Sch2}, or simply mimicking the techniques from Hopf algebra theory).\\

\noindent The following proposition is then not surprising.

\begin{Prop}\label{PropWHA} Fix $s,t\in\mathbb{R}$. Then the $U_q(\mu,\nu)$ and $\Delta_{\mu\nu}^{\upsilon}$, with $\mu,\nu,\upsilon\in\{s,t\}$, form a co-linking weak Hopf $^*$-algebra.\end{Prop}

\begin{proof} Denote by $H$ the direct sum of the $U_q(\mu,\nu)$, and by $\Delta_H$ the map as in Definition \ref{DefWHo}. Then for $x\in U_q(\mu,\nu)$, we have, by Lemma \ref{LemCoas}, \begin{eqnarray*} (\Delta_H\otimes \iota)\Delta_H(x) &=& \sum_{\upsilon,\omega} (\Delta_{\mu\upsilon}^{\omega}\otimes \iota)\Delta_{\mu\nu}^\upsilon(x) \\ &=& \sum_{\upsilon,\omega} (\iota\otimes \Delta_{\omega\nu}^{\upsilon})\Delta_{\mu\nu}^{\omega}(x) \\ &=& (\iota\otimes \Delta_H)\Delta_H(x).\end{eqnarray*}

\noindent Define now \[\varepsilon_H: H\rightarrow \mathbb{C}: U_q(\mu,\nu)\ni x_{\mu\nu} \rightarrow \delta_{\mu,\nu} \varepsilon_{\mu}(x_{\mu\nu}).\] Then it is immediately seen to form a counit for the coalgebra $(H,\Delta_H)$, since the $U_q(\mu,\nu)$ are $U_q(\mu)$-$U_q(\nu)$-bi-comodules. The `monoidality' condition A.6 in Definition 2.1 of \cite{Boh1} is also easily checked: for example, for $x\in U_q(\mu,\nu)$, $y\in U_q(\upsilon,\omega)$ and $z\in U_q(\varpi,\vartheta)$, we have \begin{eqnarray*}\varepsilon_H(xy_{(1)})\varepsilon_H(y_{(2)}z) &=& \sum_{\kappa} \varepsilon_H(xy_{(1)\upsilon\kappa})\varepsilon_H(y_{(2)\kappa\omega}z) \\ &=& \delta_{\mu,\upsilon}\delta_{\nu,\varpi} \delta_{\omega,\vartheta} \delta_{\mu,\nu}\delta_{\varpi,\vartheta} \varepsilon_{\mu}(xy_{(1)\mu\mu}) \varepsilon_{\mu}(y_{(2)\mu\mu}z) \\ &=& \delta_{\mu,\upsilon}\delta_{\nu,\varpi} \delta_{\omega,\vartheta} \delta_{\mu,\nu}\delta_{\varpi,\vartheta} \varepsilon_{\mu}(xyz) \\ &=& \varepsilon_H(xyz).\end{eqnarray*}

\noindent Similarly, it is immediately verified that $1_H = \sum_{\mu,\nu} 1_{\mu\nu}$ satisfies the `comonoidality' condition A.7 of that definition: for example \begin{eqnarray*} (\Delta_H(1_H)\otimes 1_H)(1_H\otimes \Delta_H(1_H)) &=& \sum_{\textrm{all variables}} (1_{\mu\upsilon}\otimes 1_{\upsilon\nu}\otimes 1_{\omega,\varpi})(1_{\mu'\nu'} \otimes 1_{\upsilon'\omega'}\otimes 1_{\omega'\varpi'})\\ &=& \sum_{\textrm{all variables}} 1_{\mu\upsilon} \otimes 1_{\upsilon\omega}\otimes 1_{\omega\nu} \\ &=& \Delta^{(2)}(1_H).\end{eqnarray*}

\noindent We have shown now that $(H,\Delta_H)$ is a weak bi-$^*$-algebra.\\

\noindent Finally, let us define the map \[S_H: H\rightarrow H: U_q(\mu,\nu) \ni x\rightarrow S_{\mu\nu}(x)\in U_q(\nu,\mu)\subseteq H.\] Then using Lemma \ref{LemAnt}, we get for $x\in U_q(\mu,\nu)$ that \begin{eqnarray*} x_{(1)}S_H(x_{(2)}) &=&  \sum_{\omega} x_{(1)\mu\omega}S_{\omega\nu}(x_{(2)\omega\nu}) \\ &=& \delta_{\mu,\nu} \varepsilon_{\mu}(x) \sum_{\omega} 1_{\mu\omega} \\ &=& \sum_{\omega,\varpi,\vartheta} \varepsilon_H(1_{\omega\varpi}x)1_{\varpi\vartheta} \\ &=& (\varepsilon_H\otimes \iota)(\Delta_H(1_H)(x\otimes 1)),\end{eqnarray*} proving the first antipode condition of A.8 of Definition 2.1 in \cite{Boh1}. For the condition A.9, we compute for $x\in U_q(\mu,\nu)$ that \begin{eqnarray*} S_H(x_{(1)})x_{(2)}S_H(x_{(3)}) &=& \sum_{\upsilon,\omega} S_{\mu\upsilon}(x_{(1)\mu\upsilon})x_{(2)\upsilon\omega} S_{\omega\nu}(x_{(3)\omega\nu}) \\&=& S_{\mu\nu}(x_{(1)\mu\nu})x_{(2)\nu\mu} S_{\mu\nu}(x_{(3)\mu\nu}) \\ &=& \varepsilon_{\mu}(x_{(1)\mu\mu}) S_{\mu\nu}(x_{(2)\mu\nu}) \\ &=& S_H(x).\end{eqnarray*}

\end{proof}

\noindent We now introduce the Miyashita-Ulbrich action associated to bi-$^*$-Galois objects and co-linking weak Hopf $^*$-algebras.

\begin{Def} Let $(H_{ij},\Delta_{ij}^k)$ be a co-linking weak Hopf $^*$-algebra. The \emph{left (resp.~ right) Miyashita-Ulbrich action} of the Hopf $^*$-algebra $(H_0,\Delta_0)$ (resp.~ Hopf $^*$-algebra $(H_1,\Delta_1)$) on $H_{01}$ is the left $H_0$-module (resp.~ right $H_1$-module) $^*$-algebra structure on $H_{01}$ determined as \[x\rhd y = x_{(1)01}\cdot y \cdot S_{10}(x_{(2)10}),\qquad x\in H_0,y\in H_{01},\]\[y \lhd x = S_{10}(x_{(1)10})\cdot y \cdot x_{(2)01},\qquad x\in H_1,y\in H_{01}.\]
\end{Def}

\noindent The fact that this is a left module $^*$-algebra structure means that $\rhd$ determines a left module structure of $H_0$ on $H_{01}$, which interacts with the $^*$-algebra structure of $H_{01}$ in the following way: \[ x\rhd (yz) = (x_{(1)}\rhd y)\cdot (x_{(2)}\rhd z),\]\[x\rhd (y^*) = ((S(x)^*)\rhd y)^*.\] That the above Miyashita-Ulbrich action satisfies these conditions is easily derived using the properties of the antipode $S$ of a weak Hopf $^*$-algebra, namely that it is an anti-homomorphism, satisfying $S(S(x^*)^*)=x$ for all $x$ in the weak Hopf $^*$-algebra.\\

\noindent Using the properties of the antipode, it is also not difficult to see that the above definition coincides with the usual definition of the Miyashita-Ulbrich action for bi-Galois objects (see e.g.~ \cite{Sch3}, Definition 2.1.8).\\

\noindent \emph{Note:} In the following, we will want to work with multiple bi-$^*$-Galois objects at the same time. The notion of a co-linking weak Hopf $^*$-algebra is well-adapted to this: when we have a bi-$^*$-Galois object between two Hopf $^*$-algebras $H_{-1}$ and $H_0$ and a bi-Galois object between $H_0$ and yet another Hopf $^*$-algebra $H_1$, we can group them all together, in an essentially unique way, into a `3$\times$3-co-linking weak Hopf $^*$-algebra' consisting of 9 $^*$-algebras and 27 comultiplications, or one large weak Hopf $^*$-algebra with the $H_i$ at its `corners'. Using the same techniques as above, it is then easy to see that the $U_q(\mu,\nu)$ and $\Delta_{\mu\nu}^{\upsilon}$, with indices now ranging over $\{-1,0,1\} \equiv \{-,0,+\}$, form such a 3$\times$3 co-linking weak Hopf $^*$-algebra.

\section{On the quantization of the infinitesimal action of $SU(2)$, $E(2)$ and $SU(1,1)$ on their homogeneous spaces.}

\noindent We now define Casimir elements inside our 3$\times$3 co-linking weak Hopf $^*$-algebra. We keep notation as in the previous section.

\begin{Def}\label{DefCas} For $\mu,\nu\in \{-,0,+\}$, we define the \emph{Casimir element} of $U_q(\mu,\nu)$ to be the element \begin{eqnarray*} C_{\mu\nu} &:=& EF + \lambda^2(q^{-1}\mu K^2+q \nu K^{-2})\\ &=& FE+\lambda^2(q\mu K^2+q^{-1}\nu K^{-2}).\end{eqnarray*}\end{Def}

\noindent \emph{Warning:} One should check the above equality by using the commutation relations, \emph{not} by applying the antipode. Indeed, in bi-Galois objects, the antipode is \emph{external}, and one easily checks that its naive application would violate the above equality.\\

\noindent As in the quantized enveloping algebra case, we have the following easy lemma.

\begin{Lem}\label{LemCas} The Casimir element $C_{\mu\nu}$ is a self-adjoint element in the center of $U_q(\mu,\nu)$. \end{Lem}

\begin{proof} The fact that $C_{\mu\nu}$ is self-adjoint is apparent on sight. To see that it lies in the center, we then only have to see if it commutes with $K$ and $E$. Commutation with $K$ is immediate, while for $E$, we have \begin{eqnarray*} C_{\mu\nu}E &=& EFE + \lambda^2q^{-1}\mu K^2E+ \lambda^2q \nu K^{-2}E \\ &=& E(EF-\lambda \mu K^2 + \lambda \nu K^{-2}) + \lambda^2q\mu EK^2+ \lambda^2q^{-1} \nu EK^{-2} \\ &=& E(EF -\lambda^2 (\lambda^{-1} - q)\mu K^2 + \lambda^2 (\lambda^{-1}+q^{-1})\nu K^{-2}) \\ &=& E(EF + \lambda^2 q^{-1} \mu K^2 + \lambda^2 q \nu K^{-2}) \\ &=& EC_{\mu\nu}.\end{eqnarray*}
\end{proof}

\begin{Def} For $\tau\in \mathbb{R}$, we define the $^*$-algebra $A_{\mu\nu}^{\tau}$ as the quotient $^*$-algebra of $U_q(\mu,\nu)$, obtained by evaluating $C_{\mu\nu}$ at $\tau q^{-1}\lambda^2$: \[A_{\mu\nu}^{\tau} = U_q(\mu,\nu)/(C_{\mu\nu}-\tau q^{-1}\lambda^2 1_{\mu\nu}).\]\end{Def}

\noindent The quotient will of course inherit a $^*$-structure as $C_{\mu\nu}$ is self-adjoint and $\tau$ is real.

\begin{Lem}\label{LemProj} Denote by $\pi_{\tau}$ the projection map $U_q(\mu,\nu)\rightarrow A_{\mu\nu}^{\tau}$. Then $A_{\mu\nu}^{\tau}$ inherits a left $U_q(\mu)$-module $^*$-algebra structure from $U_q(\mu,\nu)$, uniquely determined by \[x\rhd \pi_{\tau}(y) = \pi_{\tau}(x\rhd y), \qquad x\in U_q(\mu), y\in U_q(\mu,\nu).\]
\end{Lem}

\begin{proof} The Miyashita-Ulbrich action descends to \emph{any} quotient, since, by its definition, any 2-sided ideal is preserved by it.
\end{proof}

\noindent We now introduce a particular sub-$^*$-algebra of $A_{\mu\nu}^{\tau}$.\\

\begin{Not} We denote by $B_{\mu\nu}^{\tau}$ the sub-$^*$-algebra of $A_{\mu\nu}^{\tau}$ generated by the images of the elements $K_{\mu\nu}^2$ and $K_{\mu\nu}E_{\mu\nu}$ under the quotient map $\pi_{\tau}$ from $U_q(\mu,\nu)$. We denote by $X$ and $Z$ the following elements in this algebra:\begin{eqnarray*} X &=& q^{1/2}(q^{-1}-q)\,\pi_{\tau}(F_{\mu\nu}K_{\mu\nu}),\\  Z &=& \pi_{\tau}(K_{\mu\nu}^2).\end{eqnarray*} \end{Not}

\begin{Prop}\label{PropUniv} The $^*$-algebra $B_{\mu\nu}^{\tau}$ is isomorphic to the universal unital $^*$-algebra $\mathscr{B}_{\mu\nu}^{\tau}$ generated by two elements $x$ and $z$, satisfying the relations \[ z^*=z, \;\;\;\;\;\; xz = q^{2} zx, \;\;\;\;\;\; x^{*}z = q^{-2}zx^{*}\] and \begin{eqnarray*} x^* x &=& -q^2\nu+\quad \!\tau z -\quad\mu z^2, \\ xx^* &=& -q^2\nu + q^2 \tau z - q^4\mu z^2.\end{eqnarray*} An isomorphism is provided by $^*$-homomorphically extending the assignment $x\rightarrow X$ and $z\rightarrow Z$.

\end{Prop}

\begin{proof} We first verify that the elements $X$ and $Z$ of $B_{\mu\nu}^{\tau}$ satisfy the same commutation relations as $x$ and $z$. The $q$-commutation relations are of course immediate. We compute the identity for $X^*X$, and leave the other one to the reader.\\

\noindent We have \begin{eqnarray*} X^*X &=& q\lambda^{-2}\pi_{\tau}(K_{\mu\nu}E_{\mu\nu}F_{\mu\nu}K_{\mu\nu}) \\ &=& q\lambda^{-2}\pi_{\tau}(q^{-1}\lambda^2\tau-\lambda^2(q^{-1}\mu K_{\mu\nu}^2+q \nu K_{\mu\nu}^{-2}))\pi_{\tau}(K_{\mu\nu}^2) \\ &=& -q^2\nu + \tau Z -\mu Z^2.\end{eqnarray*}\\

\noindent We now prove universality. Denote by $V$ the vector space with basis vectors $e_{nm}$, $n,m\in \mathbb{Z}$, and define linear maps $\widetilde{x},\widetilde{y}$ and $\widetilde{w}$ on $V$ by the formulas \[\begin{array}{llll} \widetilde{x}\cdot e_{nm} &=& \qquad \;\, e_{n+1,m} & n\geq 0,\\ \widetilde{x}\cdot e_{nm} &=& -q^2\nu \,e_{n+1,m} + \tau q^{-2n}  \quad \;\;\, e_{n+1,m+2} - \mu q^{-4n} \quad \;\;\, e_{n+1,m+4} & n<0,\\  \widetilde{y}\cdot e_{nm} &=& -q^2\nu \,e_{n-1,m} + \tau q^{-2(n-1)} \,e_{n-1,m+2} - \mu q^{-4(n-1)}\, e_{n-1,m+4} & n>0,\\ \widetilde{y}\cdot e_{nm} &=& \qquad\;\, e_{n-1,m} & n\leq 0,\\ \widetilde{w} \cdot e_{nm} &=& \;\;\,q^{-n}\, e_{n,m+1}&.\end{array}\] Then $\widetilde{w}$ is invertible, and the set $\{\widetilde{x}^k\widetilde{w}^m,(\widetilde{y})^{k'}\widetilde{w}^{m'}\mid k,k' \in \mathbb{N}, m,m'\in \mathbb{Z} \}$ will be linearly independent, since $\widetilde{x}^k\widetilde{w}^m e_{00} = e_{km}$ and $(\widetilde{y})^{k}\widetilde{w}^{m} e_{00} = e_{-k,m}$ for $k\in \mathbb{N}$ and $m\in \mathbb{Z}$.\\

\noindent A straightforward computation shows that we can then define a unital homomorphism from $A_{\mu\nu}^{\tau}$ to $\textrm{End}(V)$ such that $X\rightarrow \widetilde{x},\pi_{\tau}(K_{\mu\nu})\rightarrow \widetilde{w}^2$ and $X^*\rightarrow \widetilde{y}$. This shows that the set $\{X^kZ^m,(X^*)^{k}Z^{m}\mid k,m \in \mathbb{N}\}$ is linearly independent in $B_{\mu\nu}^{\tau}$. Since it is easily seen that $\{x^kz^m,(x^*)^{k}z^{m}\mid k,m \in \mathbb{N}\}$ spans $\mathscr{B}_{\mu\nu}^{\tau}$ as a vector space, it follows that the natural unital $^*$-homomorphism $\mathscr{B}_{\mu\nu}^{\tau}\rightarrow B_{\mu\nu}^{\tau}$ is bijective.

\end{proof}

\noindent \emph{Remark:} Suppose that $\mu\nu\in \{-,0,+\}$ and $\tau\in \mathbb{R}$, and suppose that we are in one of the following situations: \begin{itemize} \item $\nu = -$ ; \item $\nu= 0,\tau\neq 0$; \item $\nu= \tau=0$, $\mu=-$;\item $\nu=+$, $\mu= -$;\item $\nu=+$, $\mu=0$, $\tau\neq 0$.\end{itemize} Then the $^*$-algebra $B_{\mu\nu}^{\tau}$ has a `topological implementation', in that there exists a Hilbert space $\mathscr{H}$ and a dense subset $\mathscr{D}\subseteq \mathscr{H}$ such that $B_{\mu\nu}^{\tau}$ can be represented faithfully by adjointable operators on $\mathscr{D}$, in such a way that the $^*$-operation coincides with restricting the adjoint to $\mathscr{D}$. In the remaining cases, this is impossible. Although not so difficult, we do not give a proof of this statement here, since we want to keep this paper at the level of elementary algebra. The detailed study of these spaces on a C$^*$-algebraic level, along with their further structure (see Proposition \ref{coac}), will be treated elsewhere (in as far as it has not been treated in the literature yet).\\

\begin{Prop}\label{PropInf} The left $U_q(\mu)$-module structure on $A_{\mu\nu}^{\tau}$ restricts to a left $U_q(\mu)$-module structure on $B_{\mu\nu}^{\tau}$, which hence becomes a left $U_q(\mu)$-module $^*$-algebra.

\end{Prop}

\begin{proof}

\noindent Since $A_{\mu\nu}^{\tau}$ is a left module $^*$-algebra, it is sufficient to see if the action of any of the generators $E_{\mu},F_{\mu}$ and $K_{\mu}^{\pm}$ of $U_q(\mu)$ on the elements $X,Z$ of $B_{\mu\nu}^{\tau}$ gives an element in $B_{\mu\nu}^{\tau}$. \\

\noindent Now by definition, \begin{eqnarray*} E_{\mu}\rhd X &=& q^{1/2}(q^{-1}-q) \pi_{\tau}(E_{\mu}\rhd F_{\mu\nu}K_{\mu\nu})\\ &=& q^{1/2}(q^{-1}-q)\pi_{\tau}( E_{\mu\nu}(F_{\mu\nu}K_{\mu\nu})K_{\mu\nu}^{-1} + K_{\mu\nu}^{-1}(F_{\mu\nu}K_{\mu\nu})(-qE_{\mu\nu}))\\ &=& q^{1/2}(q^{-1}-q)\pi_{\tau}( E_{\mu\nu}F_{\mu\nu} -q^2 F_{\mu\nu}E_{\mu\nu}) \\ &=& q^{1/2}(q^{-1}-q)\pi_{\tau}(q^{-1}\lambda^2\tau-\lambda^2(q^{-1}\mu K_{\mu\nu}^2+q\nu K_{\mu\nu}^{-2}) \\ && \qquad \quad -q\lambda^2\tau +q^2\lambda^2(q\mu K_{\mu\nu}^2+q^{-1}\nu K_{\mu\nu}^{-2})) \\ &=& q^{1/2}(q^{-1}-q)((q^{-1}-q)\lambda^2\tau - \lambda^2\mu(q^{-1}-q^3)Z) \\ &=& q^{1/2}\tau - q^{1/2} (1+q^2)\mu Z,\\
\\
F_{\mu}\rhd X &=& q^{1/2}(q^{-1}-q) \pi_{\tau}(F_{\mu}\rhd F_{\mu\nu}K_{\mu\nu})\\ &=& q^{1/2}(q^{-1}-q)\pi_{\tau}( F_{\mu\nu}(F_{\mu\nu}K_{\mu\nu})K_{\mu\nu}^{-1} + K_{\mu\nu}^{-1}(F_{\mu\nu}K_{\mu\nu})(-q^{-1}F_{\mu\nu}))\\ &=& 0, \end{eqnarray*}
\begin{eqnarray*}K_{\mu}\rhd X &=& q^{1/2}(q^{-1}-q) \pi_{\tau}(K_{\mu}\rhd F_{\mu\nu}K_{\mu\nu})\\ &=& q^{1/2}(q^{-1}-q)\pi_{\tau}( K_{\mu\nu}(F_{\mu\nu}K_{\mu\nu})K_{\mu\nu}^{-1}) \\ &=& q^{-1}X,\\ \\
E_{\mu}\rhd Z &=& \pi_{\tau}(E_{\mu}\rhd K_{\mu\nu}^2) \\ &=& \pi_{\tau}(E_{\mu\nu}(K_{\mu\nu}^2)K_{\mu\nu}^{-1} + K_{\mu\nu}^{-1}(K_{\mu\nu}^2)(-qE_{\mu\nu})) \\ &=& (q^{-1}-q)\pi_{\tau}(K_{\mu\nu}E_{\mu\nu})  \\ &=& q^{-1/2} X^*, \\ K_{\mu}\rhd Z &=& Z. \end{eqnarray*}

\noindent Using the module $^*$-structure, we also find \begin{eqnarray*} F_{\mu}\rhd Z &=& F_{\mu}\rhd Z^* \\ &=& (S_{\mu}(F_{\mu})^* \rhd Z)^* \\ &=& -q^{-1} (E_{\mu}\rhd Z)^* \\&=& -q^{-3/2} X.\end{eqnarray*}

\end{proof}

\noindent From the proof of the foregoing proposition, we get the following formulas for the $U_q(\mu)$-module structure on $B_{\mu\nu}^{\tau}$, using the module $^*$-structure for the identities in the final column:

\begin{center}\begin{tabular}{l|ccccc} $\rhd$ & $X$ && $Z$ && $X^*$ \\ \hline \\ $E_{\mu}$ & $q^{1/2}(\tau -  (1+q^2)\mu Z)$ &&  $q^{-1/2} X^*$ && $0$ \\\\ $K_{\mu}$ & $q^{-1} X$ && $Z$ && $qX^*$ \\\\$F_{\mu}$ & $0$ && $-q^{-3/2} X$ &&$-q^{-1/2}(\tau - (1+q^2)\mu Z) $\end{tabular}\end{center}\vspace{0.4cm}

\noindent We now want to show that this module structure is the `infinitesimal' module structure associated to some coaction by a dual Hopf $^*$-algebra. We first introduce the relevant objects.

\begin{Def} The \emph{algebra of polynomial functions on $SL_q(2,\mathbb{C})$} (\cite{Kli1}) is defined as the unital algebra $\textrm{Pol}(SL_q(2,\mathbb{C}))$ generated by four generators $a,b,c,d$ satisfying the relations \[ ab = qba, \quad ac = qca, \quad bd = qdb, \quad cd = qdc, \quad bc = cb\] and \[da -q^{-1}cb = 1, \quad ad - q bc = 1.\]

\noindent It can be made into a Hopf algebra by endowing it with the comultiplication map $\Delta$ satisfying \[\left\{\begin{array}{l} \Delta(a) = a\otimes a + c \otimes b \\ \Delta(b) = b\otimes a+ d\otimes b, \\  \Delta(c) = a\otimes c+ c\otimes d, \\
\Delta(d) = b\otimes c + d\otimes d. \end{array}\right.\]

\noindent The algebra $\textrm{Pol}(SL_q(2))$ can be endowed with the $^*$-structure $a^*=d, b^* = -q^{-1}c$, in which case it becomes a Hopf $^*$-algebra which we denote as $(\textrm{Pol}(SU_q(2)),\Delta_+) = \textrm{Pol}_q(+,+)$ (\cite{Wor1}), and whose generators we denote then as $a_+$ and $b_+$.\\

\noindent The algebra $\textrm{Pol}(SL_q(2))$ can also be endowed with the $^*$-structure $a^*=d,b^*=q^{-1}c$, in which case it becomes a Hopf $^*$-algebra which we denote as $(\textrm{Pol}(SU_q(1,1)),\Delta_-)=\textrm{Pol}_q(-,-)$ (\cite{Kor2}), and whose generators we denote then as $a_-$ and $b_-$.\\

\vspace{1cm}

\noindent The $^*$-algebra of \emph{polynomial functions on the quantum $\widetilde{E}(2)$ group} (\cite{Vak1},\cite{Cel1},\cite{Koe1}) is defined as the universal unital $^*$-algebra \[\textrm{Pol}_q(0,0)=\textrm{Pol}(\widetilde{E}_q(2))\] generated by elements $a_0$ and $b_0$, subject to the relations
\[\left\{\begin{array}{lllclll} a_0^*a_0 = 1 && \!\!\!\!\!\!\!\!\!\!\!\!a_0b_0 = q\;\;\;\,b_0a_0 \\a_0a_0^* = 1 &&\!\!\!\!\!\!\!\!\!\!\!\! a_0^*b_0 = q^{-1}b_0a_0^* \\ &\!\!\!\!\!\!\!\!b_0b_0^* = b_0^*b_0.\end{array}\right.\] We can make it into a Hopf $^*$-algebra  by endowing it with the comultiplication map $\Delta_0$ satisfying \begin{eqnarray*} \Delta_0(a_0) &=& a_0\otimes a_0 \\ \Delta_0(b_0) &=& b_0\otimes a_0 + a_0^*\otimes b_0.\end{eqnarray*}

\end{Def}

\noindent We recall now that a non-degenerate pairing between Hopf $^*$-algebras $(H,\Delta_H)$ and $(K,\Delta_K)$ (\cite{VDae2}) consists of a non-degenerate bilinear map \[\langle\;\cdot\;,\;\cdot\;\rangle: H\times K\rightarrow \mathbb{C}\] such that \[ \langle \Delta_H(x),y\otimes z \rangle = \langle x,yz\rangle \qquad \langle x\otimes y,\Delta_K(z) \rangle = \langle xy,z\rangle,\] and such that \[\langle x^*,y\rangle = \overline{\langle x,S(y)^*\rangle}, \qquad \langle x,y^*\rangle = \overline{\langle S(x)^*,y\rangle}.\]

\begin{Prop}\label{PropPair} For $\mu\in \{-,0,+\}$, there is a non-degenerate pairing between the Hopf $^*$-algebras $U_q(\mu)$ and $\textrm{Pol}_q(\mu)$, uniquely determined by the formulas \[ \langle K,a \rangle = q^{-1/2}, \quad  \langle K,a^*\rangle = q^{1/2}, \quad \langle E,b \rangle = 1, \quad \langle F,(-qb^*)\rangle = 1,\] while all other possible pairings between generators are assigned zero.\end{Prop}

\begin{proof} For $\mu=+$ and $\mu=-$, this is well-known (see \cite{Kli1}, Theorem 21 and the discussion following it. We remark that in the case $\mu=-$, their generator $F$ corresponds to our $-F_-$). For $\mu=0$, it follows from \cite{Koe1}, Corollary 3.5.\\

\noindent The Proposition can also be checked directly, except for the non-degeneracy statement, by applying the argument in the beginning of section 4 of \cite{VDae2}.
\end{proof}

\begin{Prop}\label{coac} There exists a right coaction $\gamma: B_{\mu\nu}^{\tau}\rightarrow B_{\mu\nu}^{\tau}\otimes \textrm{Pol}_q(\mu)$ such that \[\begin{array}{lcrcrcrcr} \gamma(X) &=& -q \mu X^*\otimes b_{\mu}^2 &-& q(1+q^2)\mu Z \otimes b_{\mu}a_{\mu} &+&X\otimes a_{\mu}^2  &+& q\tau 1\otimes b_{\mu}a_{\mu},\\
\gamma(Z) &=&  X^* \otimes a_{\mu}^*b_{\mu} &+& Z\otimes (1-(1+q^2)\mu b_{\mu}^*b_{\mu}) &+& X\otimes b_{\mu}^*a_{\mu} &+& \tau 1\otimes b_{\mu}^*b_{\mu},\\
\gamma(X^*) &=& X^*\otimes (a_{\mu}^*)^2  &-& q(1+q^2)\mu Z \otimes a_{\mu}^*b_{\mu}^*  &-&q \mu X\otimes (b_{\mu}^*)^2 &+& q \tau 1\otimes a_{\mu}^*b_{\mu}^*,\end{array}\] and such that, for any $x\in U_q(\mu)$ and $y\in B_{\mu\nu}^{\tau}$, \[x\rhd y = (\iota\otimes \langle \;\cdot\;,x\rangle)\gamma(y).\]

\end{Prop}

\begin{proof} Put \begin{eqnarray*} \omega_1&=&-q^{-1}X,\\ \omega_0&=& (1+q^2)^{1/2}( \mu Z- (1+q^2)^{-1}\tau), \\ \omega_{-1} &=& \mu X^*,\end{eqnarray*} and let $V$ be the linear span of the $\omega_i$ and 1. Then by the formulas for $\gamma$ stated in the proposition, together with the formula $\gamma(1)=1\otimes 1$, we can certainly construct $\gamma$ as a linear map $\gamma: V\rightarrow V\otimes \textrm{Pol}_q(\mu)$, and it then satisfies \begin{eqnarray*} \gamma(\omega_1) &=&  \omega_{-1} \otimes b_{\mu}^2 + (1+q^2)^{1/2} \omega_0 \otimes b_{\mu}a_{\mu} + \omega_1\otimes a_{\mu}^2,\\
\gamma(\omega_0) &=&  (1+q^{2})^{1/2}  \omega_{-1} \otimes a_{\mu}^*b_{\mu} + \omega_0\otimes (1-(1+q^2)\mu b_{\mu}^*b_{\mu}) - (1+q^2)^{1/2}q\mu \omega_1\otimes b_{\mu}^*a_{\mu},\\
\gamma(\omega_{-1}) &=&  \omega_{-1}\otimes (a_{\mu}^*)^2- q\mu (1+q^2)^{1/2} \omega_0 \otimes a_{\mu}^*b_{\mu}^*  + q^2 \mu^2 \omega_{1} \otimes (b_{\mu}^*)^2 ,
\end{eqnarray*}

\noindent so

\[  \gamma \left(\begin{array}{ccc} \omega_{-1} & \omega_0 & \omega_1\end{array}\right) = \left(\begin{array}{ccc} \omega_{-1} & \omega_0 & \omega_1\end{array}\right) \otimes \left(\begin{array}{ccc} d_{\mu}^2 & (1+q^2)^{1/2} d_{\mu}b_{\mu} & b_{\mu}^2 \\ (1+q^2)^{1/2} d_{\mu}c_{\mu} & 1+(q+q^{-1}) b_{\mu}c_{\mu} & (1+q^2)^{1/2} b_{\mu}a_{\mu} \\ c_{\mu}^2 & (1+q^2)^{1/2} c_{\mu}a_{\mu} & a_{\mu}^2  \end{array}\right),\] where $c_{\mu} = -q \mu b_{\mu}^*$ and $d_{\mu} = a_{\mu}^*$. But  $\left(\begin{array}{ccc} d^2 & (1+q^2)^{1/2} db & b^2 \\ (1+q^2)^{1/2} dc & 1+(q+q^{-1}) bc & (1+q^2)^{1/2} ba \\ c^2 & (1+q^2)^{1/2} ca & a^2  \end{array}\right)$ is a spin 1 corepresentation of $\textrm{Pol}(SL_q(2,\mathbb{C}))$ (see e.g.~ the concrete form used in \cite{Mas2}). Since \[ \textrm{Pol}(SL_q(2,\mathbb{C}))\rightarrow \textrm{Pol}_q(\mu): x\rightarrow x_{\mu},\quad x\in\{a,b,c,d\}\] extends to a homomorphism of Hopf algebras, we have that $(V,\gamma)$ will be a right $\textrm{Pol}_q(\mu)$-comodule.\\

\noindent Since the only pairings between elements of $\{E_{\mu},F_{\mu},K_{\mu}\}$ and elements of \[\{1,b_{\mu}^2,b_{\mu}a_{\mu},a_{\mu}^2,a_{\mu}^*b_{\mu},b_{\mu}^*b_{\mu},b_{\mu}^*a_{\mu},(a_{\mu}^*)^2,a_{\mu}^*b_{\mu}^*,(b_{\mu}^*)^2\}\] which are \emph{not} zero are \[\langle E_{\mu}, b_{\mu}a_{\mu} \rangle = q^{-1/2}= \langle E_{\mu},a_{\mu}^*b_{\mu}\rangle,\quad  \langle F_{\mu}, b_{\mu}^*a_{\mu} \rangle = -q^{-3/2} = \langle F_{\mu},a_{\mu}^*b_{\mu}^*\rangle\] and \[\langle K_{\mu},a_{\mu}^2 \rangle = q^{-1},\quad \langle K_{\mu},1\rangle = 1, \quad \langle K_{\mu},(a_{\mu}^*)^2 \rangle = q,\] it is quite immediately verified that for any $x\in \{E_{\mu},F_{\mu},K_{\mu}\}$ and $y\in V$, we have \[x\rhd y = (\iota\otimes \langle \;\cdot\;,x\rangle)\gamma(y).\] Then this formula is of course true for any $x\in U_q(\mu)$.\\

\noindent But since $B_{\mu\nu}^{\tau}$ is a left module $^*$-algebra for $U_q(\mu)$, and since the pairing between $\textrm{Pol}_q(\mu)$ and $U_q(\mu)$ is non-degenerate, we get that $\gamma$ can be extended to a right coaction on the $^*$-algebra generated by $V$. In case $\mu\neq 0$, this then proves the proposition. In case $\mu=0$ but $\tau\neq 0$, we have $Z=\tau^{-1}(X^*X +q^2\nu)$, and applying $\gamma$ to this we again get the formula for $\gamma(Z)$ in the proposition.\\

\noindent Finally, in case $\mu=\tau=0$, one verifies directly that $\gamma$ as given in the proposition, restricted to the linear span of $X^*,Z$ and $X$, gives a $\textrm{Pol}_q(0)$-comodule structure, whose associated $U_q(0)$-module structure coincides with the one coming from $B_{0\nu}^0$. The same argument as before then lets us conclude that also in this case, $\gamma$ can indeed be extended to a coaction.

\end{proof}

\noindent \emph{Remarks}:\begin{enumerate}\item The case $\mu=0$ in the previous Proposition could also have been derived from the $\mu\neq 0$-case by allowing $\mu$ to take values in $\mathbb{R}$ again and then letting $\mu\rightarrow 0$.
\item In case $\mu=+$, we obtain in this way the well-known \emph{Podle\'{s} quantum spheres} (\cite{Pod1}), with their natural coaction by $\textrm{Pol}(SU_q(2))$.
\end{enumerate}

\begin{Prop} For all $\mu,\nu \in\{-,0,+\}$, and $\tau\in \mathbb{R}$, the coaction $\gamma$ on $B_{\mu\nu}^{\tau}$ is ergodic: if $x\in B_{\mu\nu}^{\tau}$ satisfies $\gamma(x) = x\otimes 1$, then $x \in \mathbb{C}1$.\end{Prop}
\begin{proof}

\noindent Let $V$ be the space of fixed elements for $\gamma$: $V$ consists of the elements $x\in B_{\mu\nu}^{\tau}$ such that $\gamma(x) = x\otimes 1$. Then also $y\rhd x = \varepsilon_{\mu}(y) x$ for all $y\in U_q(\mu)$.\\

\noindent Now from the proof of Proposition \ref{PropUniv}, we know that $B_{\mu\nu}^{\tau}$ has a basis consisting of vectors of the form $X^nZ^m$ and $(X^*)^{n}Z^m$ with $n,m\in \mathbb{N}$. Since $K\rhd (X^n Z^m) = q^{-n} X^nZ^m$ and $K\rhd (X^*)^n Z^m = q^n(X^*)^nZ^m$, we see that necessarily any element in $V$ must be a polynomial in $Z$. But since $\gamma(Z) = X^*\otimes A+Z\otimes B + X\otimes b_{\mu}^*a_{\mu}+1\otimes  D$ for certain elements $A,B,D\in \textrm{Pol}_q(\mu)$, we see that if $P(Z)$ is a polynomial in $Z$ of degree $k\geq 0$, and $\omega$ is the functional on $B_{\mu\nu}^\tau$ such that $\omega(X^nZ^m) = \delta_{n,k}\delta_{m,0}$ and $\omega((X^*)^nZ^m)=0$, we get $(\omega\otimes \iota)\gamma(P(Z)) \neq 0$. So if $P(Z)\in V$, necessarily $k\leq0$. Hence $\gamma$ is ergodic.

\end{proof}

\noindent By this proposition, one may look upon the $B_{\mu\nu}^{\tau}$ formally as well-behaving function spaces on quantum homogeneous spaces for the respective quantum group associated to $\textrm{Pol}_q(\mu)$.\\

\noindent Let us now determine which of the above coactions arise as coideals. The case $\mu=+$ is of course well-known, but we are not sure if the other cases have been dealt with explicitly in the literature. \\

\begin{Prop} \begin{enumerate}\item For $\tau\neq 0$, $B_{\mu\nu}^{\tau}$ is isomorphic to a right $^*$-coideal of $\textrm{Pol}_q(\mu)$ iff $\nu \leq 0$.
\item For $\tau = 0$, $B_{\mu\nu}^{\tau}$ is isomorphic to a right $^*$-coideal of $\textrm{Pol}_q(\mu)$ iff $\nu = -$.
\end{enumerate}
\end{Prop}

\begin{proof}

\noindent Let us first treat the case $\mu \neq 0$.\\

\noindent Let $r,s,t$ be complex numbers and put \[  \left(\begin{array}{ccc} \omega_{-1} & \omega_0 & \omega_1\end{array}\right) := \left(\begin{array}{ccc} r & s & t \end{array}\right) \cdot \left(\begin{array}{ccc} d_{\mu}^2 & (1+q^2)^{1/2} d_{\mu}b_{\mu} & b_{\mu}^2 \\ (1+q^2)^{1/2} d_{\mu}c_{\mu} & 1+(q+q^{-1}) b_{\mu}c_{\mu} & (1+q^2)^{1/2} b_{\mu}a_{\mu} \\ c_{\mu}^2 & (1+q^2)^{1/2} c_{\mu}a_{\mu} & a_{\mu}^2  \end{array}\right),\] where again $c_{\mu} = -q \mu b_{\mu}^*$ and $d_{\mu} = a_{\mu}^*$. Then one can compute that $\tilde{X} = -q\omega_1$, $\tilde{Y}= \mu \omega_{-1}$ and $\tilde{Z} = \mu (1+q^2)^{-1/2}(\omega_0-s)$ satisfy the commutation relations \[\widetilde{X}\widetilde{Z} = q^2 \widetilde{Z}\widetilde{X}\qquad \widetilde{Y}\widetilde{Z} = q^{-2}\widetilde{Z}\widetilde{Y},\] and \begin{eqnarray*} \tilde{Y}\tilde{X} &=& -q\mu rt - \;\;\;\, (1+q^2)^{1/2}s \tilde{Z} - \;\;\;\mu  \tilde{Z}^2,\\ \tilde{X}\tilde{Y} &=& -q\mu rt - q^2(1+q^2)^{1/2}s \tilde{Z} - q^{4}\mu \tilde{Z}^2,\end{eqnarray*} and that the unital algebra generated by $\tilde{Z}$, $\tilde{X}$ and $\tilde{Y}$ is a right coideal, with a basis consisting of elements of the form $\tilde{X}^n\tilde{Z}^m$ and $\tilde{Y}^n\tilde{Z}^m$ whenever not all $r,s,t$ are zero. Moreover, any three-tuple of elements in $\textrm{Pol}_q(\mu)$ on which $\textrm{Pol}_q(\mu)$ coacts by the above spin 1 representation is necessarily of the form  $\left(\begin{array}{ccc} \omega_{-1} & \omega_0 & \omega_1\end{array}\right)$ for certain $r,s$ and $t$. See the first section of \cite{Mas2} for some of these statements in the setting of left coideals, and also \cite{Kli1}, section 4.5.\\

\noindent From the proof of Proposition \ref{coac}, we obtain then that if $B_{\mu\nu}^{\tau}$ is to be isomorphic to a right $^*$-coideal of $\textrm{Pol}_q(\mu)$, it is necessary and sufficient that there exist $r,s$ and $t$, at least one of which is non-zero, such that the application \begin{eqnarray*} X&\rightarrow&  \tilde{X},\\ Z&\rightarrow& \widetilde{Z} + \mu(1+q^2)^{-1/2}(s + (1+q^2)^{-1/2}\tau) \\ X^*&\rightarrow & \tilde{Y}.\end{eqnarray*} extends to a $^*$-homomorphism $B_{\mu\nu}^{\tau}\rightarrow \textrm{Pol}_q(\mu)$. Then from the fact that $\tilde{Z}$ and $\tilde{X}$ $q^2$-commute, we conclude that necessarily $s= -(1+q^2)^{-1/2}\tau$, and from the fact that $\tilde{Y}$ should be $\tilde{X}^*$, we get that $r=-q\mu \overline{t}$. But comparing the commutation relation between $\tilde{X}$ and $\tilde{X}^*$, we see that it is then necessary and sufficient that $q^2|t|^2 = -q^2\nu$, which is of course only possible if $\nu< 0$ or $\nu=0$. In the latter case, the situation $\tau=0$ should be excluded since otherwise $r=s=t=0$.\\

\noindent Let us now consider the case $\mu =0$. Suppose that there exists an equivariant $^*$-isomorphism $\pi$ from $B_{0\nu}^{\tau}$ into $\textrm{Pol}_q(0)$. By the well-known formula $(v+w)^r = \sum_{k=0}^r \qbin{r}{k} w^k v^{r-k}$ between $q^{2}$-commuting variables, we get that, for any $m\in \mathbb{Z},k,l\in \mathbb{N}$, we have \[ \Delta_0(a_0^mb_0^k(b_0^*)^l) = \sum_{r=0}^k\sum_{s=0}^l c_{r,s}^{(k,l,m)} a_0^{m-r+l-s}b_0^{k-r}(b_0^*)^s \otimes a_0^{m+k-r-s} b_0^r(b_0^*)^{l-s},\] for certain non-zero numbers $c_{r,s}^{(k,l,m)}$. From this, it is not difficult to conclude, by the equivariance of $\pi$, that $\pi(X)$ should be of the form $\theta_1 a_0^2 +\tau a_0b_0$ for some $\theta_1\in \mathbb{C}$.\\

\noindent If now further $\tau\neq 0$, then $Z = \tau^{-1}(X^*X+q^2\nu)$, and hence \[\pi(Z) = \tau b_0^*b_0 +\overline{\theta_1}a_0^*b_0 + \theta_1 b_0^*a_0 + \tau^{-1}(|\theta_1|^2 +q^2\nu).\] By the $q^2$-commutation between $\pi(X)$ and $\pi(Z)$, we find the condition $|\theta_1|^2 = -q^2\nu$. Hence $\nu \leq 0$. If $\tau = 0$, then we should have $X^*X = -q^2\nu$. But in this case $\pi(X)^*\pi(X) = |\theta_1|^2$, from which we again obtain the condition $|\theta_1|^2 = -q^2\nu$. Hence $\nu  = -$, since $\pi(X)\neq 0$.\\

\noindent Finally, we have to show that such a $\pi$ exists when $\tau\neq 0$ and $\nu \leq 0$, and when $\tau = 0$ and $\nu=-$. Suppose we are in one of these cases, and consider the following elements of $\textrm{Pol}_q(0)$: \[\widetilde{X} = - q \nu a_0^2 +\tau a_0b_0,\]\[\widetilde{Z} = -q\nu a_0^*b_0 + \tau b_0^*b_0 -  q\nu b_0^*a_0.\] Then it follows again from a straightforward computation that $\widetilde{X}$, $\widetilde{Z}$ and $\widetilde{X}^*$ satisfy the same commutation relations as $X,Z$ and $X^*$, so that we can construct $\pi: B_{\mu\nu}^{\tau} \rightarrow \textrm{Pol}_q(0)$ by $^*$-homomorphically extending the application $X\rightarrow \widetilde{X}$ and, in case $\tau =0$, $Z\rightarrow \widetilde{Z}$. Moreover, $\pi$ is also immediately verified to be equivariant (since this condition has only to be checked on generators).\\

\noindent Now the associated infinitesimal left action of $U_q(\mu)$ on $\textrm{Pol}_q(\mu)$ will satisfy $K_{\mu} \rhd \widetilde{X}^k\widetilde{Z}^l  = q^{-k}\widetilde{X}^k\widetilde{Z}^l$ and $K_{\mu} \rhd (\widetilde{X}^*)^k\widetilde{Z}^l = q^k (\widetilde{X}^*)^k\widetilde{Z}^l$, so that the spaces $\widetilde{X}^k \mathbb{C}\lbrack \widetilde{Z}\rbrack$ and $(\widetilde{X}^*)^k\mathbb{C}\lbrack \widetilde{Z}\rbrack$ will be linearly independent. Keeping track of the highest (or lowest) power of $a_0$, we can then easily see that the $\widetilde{X}^k\widetilde{Z}^l$ and $(\widetilde{X}^*)^k\widetilde{Z}^l$ are all linearly independent. From this, it follows that $\pi$ is an isomorphism, and we are done.

\end{proof}

\noindent The following proposition determines when two $B_{\mu\nu}^{\tau}$ are equivariantly isomorphic.

\begin{Prop} For $\mu,\nu,\nu'\in\{-,0,+\}$ and $\tau,\tau'\in \mathbb{R}$, we have $B_{\mu\nu}^{\tau}$ equivariantly isomorphic with $B_{\mu\nu'}^{\tau'}$ iff the following conditions hold: denoting $F_+ = F_- =\{-1,1\}$ and $F_0 = \mathbb{R}_0$, we should have \begin{itemize}\item $\nu = \nu'$, and
\item there exists $\theta  \in F_{\nu}$ with $\tau = \theta \tau'$.

\end{itemize}

\end{Prop}

\begin{proof}

\noindent Let $\pi: B_{\mu\nu}^{\tau}\rightarrow B_{\mu\nu'}^{\tau'}$ be an equivariant $^*$-isomorphism. For clarity, we now index the generators $X$ and $Z$ of $B_{\mu\nu}^{\tau}$ by the indices $\nu$ and $\tau$.\\

\noindent By the formulas following Proposition \ref{PropInf}, we easily deduce, using the module algebra structure, that $K_{\mu} \rhd X_{\nu\tau}^kZ_{\nu\tau}^l  = q^{-k}X_{\nu\tau}^kZ_{\nu\tau}^l$ and $K_{\mu} \rhd (X_{\nu\tau}^*)^kZ_{\nu\tau}^l = q^k (X_{\nu\tau}^*)^kZ_{\nu\tau}^l$ for $k,l\in \mathbb{N}$. By the proof of Proposition \ref{PropUniv}, the $X_{\nu\tau}^kZ_{\nu\tau}^l$ and $(X_{\nu\tau}^*)^kZ_{\nu\tau}^l$ form a basis for $B_{\mu\nu}^{\tau}$. So by the equivariance of $\pi$, we obtain that there must exist non-zero polynomial functions $P_X$ and $P_Z$ over $\mathbb{C}$ such that $\pi(X_{\nu\tau}) = X_{\nu'\tau'} P_X(Z_{\nu'\tau'})$ and $\pi(Z_{\nu\tau}) = P_Z(Z_{\nu'\tau'})$. Since further $\pi(X_{\nu\tau})$ and $\pi(Z_{\nu\tau})$ $q^2$-commute, and $Z_{\nu'\tau'}$ is transcendental over $\mathbb{C}$, we must have $P_Z(Z_{\nu'\tau'}) = \theta Z_{\nu'\tau'}$ for $\theta \in \mathbb{C}_0$. Since $\pi(Z_{\nu\tau})$ is self-adjoint, we must have $\theta\in \mathbb{R}_0$.\\

\noindent Now since $\pi$ is equivariant, we see by comparing the coefficients of $b_{\mu}^*a_{\mu}$ in $\gamma(Z)$, where $\gamma$ is the coaction from Proposition \ref{coac}, that also $\pi(X_{\mu\nu}) = \theta X_{\mu'\nu'}$ (so that $\pi$ is completely determined by $\theta$). Then using that $\pi(X_{\nu\tau})^*\pi(X_{\nu\tau}) = \pi(X_{\nu\tau}^*X_{\nu\tau})$, we get that \[-q^2\theta^2 \nu' + \tau' \theta^2 Z_{\nu'\tau'} -  \mu \theta^2Z_{\nu'\tau'}^2 = -q^2\nu +\tau \theta Z_{\nu'\tau'} -  \mu\theta^2 Z_{\nu'\tau'}^2 .\] From this, we arrive at the conclusion that $\nu = \nu'$ and $\tau = \theta \tau'$, and that, when $\nu\neq 0$, we further have $\theta = \pm 1$. This proves one half of the proposition.\\

\noindent Now suppose $\nu \in \{-,0,+\}$ and $\tau'\in \mathbb{R}$. Take $\theta \in F_{\nu}$, and put $\tau = \theta \tau'$. Then it is immediately seen that $\theta X_{\nu\tau'}$, $\theta Z_{\nu\tau'}$ and $\theta X_{\nu\tau'}^*$ satisfy the same commutation relations as $X_{\nu\tau}$, $Z_{\nu\tau}$ and $X_{\nu\tau}^*$, so by Proposition \ref{PropUniv}, we can extend the application \begin{eqnarray*} X_{\nu\tau}&\rightarrow& \theta X_{\nu\tau'} \\ Z_{\nu\tau} &\rightarrow& \theta Z_{\nu\tau'}\end{eqnarray*} to a $^*$-isomorphism between $B_{\mu\nu}^{\tau}$ and $B_{\mu\nu}^{\tau'}$. From the expressions in Proposition \ref{coac}, it is also immediately seen that this isomorphism is equivariant. This concludes the proof.

\end{proof}

\noindent We end with the following proposition which explains how the corresponding results for the \emph{right} Miyashita-Ulbrich action can be obtained from those for the left one.

\begin{Prop} Let $D_{\mu\nu}^{\tau}$ be the sub-$^*$-algebra generated by the images of $E_{\mu\nu}K_{\mu\nu}^{-1}$ and $K_{\mu\nu}^{-2}$ inside $A_{\mu\nu}^{\tau}$. Then the quotient right Miyashita-Ulbrich action of $U_q(\nu)$ on $A_{\mu\nu}^{\tau}$ restricts to $D_{\mu\nu}^{\tau}$, and there then exists an anti-isomorphism $\Theta:B_{\nu\mu}^{\tau}\rightarrow D_{\mu\nu}^{\tau}$, non-$^*$-preserving but satisfying $\Theta(x^*)= (\Theta^{-1}(x))^*$ for $x\in B_{\nu\mu}^{\tau}$, such that \[ \Theta(x\rhd y) = \Theta(y) \lhd S_{\nu}(x), \quad x\in U_q(\nu), y\in B_{\nu\mu}^{\tau}.\]\end{Prop}

\begin{proof} We have that $S_{\nu\mu}(C_{\nu\mu}) = C_{\mu\nu}$ for the Casimir elements of Definition \ref{DefCas}, so $S_{\nu\mu}$ descends to an anti-isomorphism $\Theta:A_{\nu\mu}^{\tau}\rightarrow A_{\mu\nu}^{\tau}$ satisfying $\Theta(x)^* = (\Theta^{-1}(x))^*$ for $x\in A_{\nu\mu}^{\tau}$. It is then immediately seen that $\Theta$ restricts to a bijection $B_{\nu\mu}^{\tau}\rightarrow D_{\mu\nu}^{\tau}$.\\

\noindent Now for $y\in U_q(\nu,\mu)$ and $x\in U_q(\nu)$, we have, using the notation from Lemma \ref{LemProj}, \begin{eqnarray*} \Theta(x\rhd \pi_{\tau}(y)) &=&  \Theta(\pi_{\tau}(x_{(1)\nu\mu}yS_{\mu\nu}(x_{(2)\mu\nu})))\\ &=&\pi_{\tau}(S_{\nu\mu}(x_{(1)\nu\mu}yS_{\mu\nu}(x_{(2)\mu\nu})))  \\ &=& \pi_{\tau}(S_{\nu\mu}((S_{\nu}(x))_{(1)\nu\mu})S_{\nu\mu}(y)(S_{\nu}(x))_{(2)\mu\nu})) \\ &=& \pi_{\tau}(S_{\nu\mu}(y))\lhd S_{\nu}(x) \\ &=& \Theta(\pi_{\tau}(y))\lhd S_{\nu}(x),\end{eqnarray*} where we have used that \[(S_{\mu\upsilon}\otimes S_{\upsilon\nu})\circ \Delta_{\mu\nu}^{\upsilon} = (\Delta_{\nu\mu}^{\upsilon})^{\textrm{op}} \circ S_{\mu\nu}\] on $U_q(\mu,\nu)$, which follows straightforwardly from the fact that the antipode on the associated co-linking weak Hopf $^*$-algebra flips the comultiplication.\\

\noindent This concludes the proof.

\end{proof}

\end{document}